\documentclass[a4paper,12pt]{article}
\usepackage{amssymb}

\newcommand{\R}{\mathbb{R}}

\newcommand{\beq}{\begin{equation} }
\newcommand{\eqq}{\end{equation} }
\newcommand{\cuad}{{\sqcap\kern-.68em\sqcup}}

\newtheorem{definition}{Definition}[section]
\newtheorem{teo}{Theorem}[section]

\newtheorem{proposition}{Proposition}[section]

\newtheorem{lemma}{Lemma}[section]

\newtheorem{remark}{Remark}[section]
\newcommand{\bremark}{\begin{remark} \em}
\newcommand{\eremark}{\end{remark} }

\def\beeq{\begin{equation}}
\def\eeq{\end{equation}}
\newcommand{\begeqaet}{\begin{eqnarray*}}
\newcommand{\eneqaet}{\end{eqnarray*}}

\newcommand{\ve}{\varepsilon}

\newcommand{\Om}{{\Omega}}

\begin{document}
\begin{center}{\bf Multiple positive solutions for nonlinear critical fractional elliptic equations involving sign-changing weight functions}\medskip

\bigskip

\bigskip

{ Alexander Quaas and Aliang Xia}

Departamento de  Matem\'atica,  Universidad T\'ecnica Federico Santa Mar\'{i}a

Casilla: V-110, Avda. Espa\~na 1680, Valpara\'{\i}so, Chile.

 {\sl ( alexander.quaas@usm.cl and  aliangxia@gmail.com)}
\end{center}

\bigskip

\begin{abstract}
In this article, we prove the  existence and multiplicity of positive solutions for the following fractional elliptic equation with sign-changing weight functions:
\begin{eqnarray*}
\left\{\begin{array}{l@{\quad }l}
(-\Delta)^\alpha u= a_\lambda(x)|u|^{q-2}u+b(x)|u|^{2^*_\alpha-1}u &{\rm in}\,\,\Omega,\\
u=0\,\,&{\rm in}\,\,\R^N\setminus\Omega,
 \end{array}
 \right.
 \end{eqnarray*}
where $0<\alpha<1$, $ \Omega $ is a bounded domain with smooth boundary in $ \R^N $ with $N>2\alpha$ and $ 2^*_\alpha=2N/(N-2\alpha)$ is the fractional critical Sobolev exponent. Our multiplicity results are based on studying
the decomposition of the Nehari manifold and  the Ljusternik-Schnirelmann category.
\end{abstract}
\date{}


\setcounter{equation}{0}
\section{ Introduction}

The fractional Laplacian has attracted much attention recently. It has applications in mathematical physics, biological modeling and mathematical finances and so on. 
Especially, it appears in turbulence and water wave, anomalous dynamics, flames propagation and chemical reactions in liquids, population dynamics, geophysical fluid dynamics, and American options in finance. For more details and applications, see \cite{A, B1, CT1, MT, V, VIKH} and references therein.

In this paper we focus our attention on critical fractional elliptic problems involving sign-changing functions. More precise, we consider the following elliptic equation involving the fractional Laplacian:
\begin{eqnarray}\label{1.1}
\left\{\begin{array}{l@{\quad }l}
(-\Delta)^\alpha u= a_\lambda(x)|u|^{q-2}u+b(x)|u|^{2^*_\alpha-1}u &{\rm in}\,\,\Omega,\\
u=0\,\,&{\rm in}\,\,\R^N\setminus\Omega,
 \end{array}
 \right.
 \end{eqnarray}
 where $ \Omega $ is a bounded domain with smooth boundary in $ \R^N $, $ 0<\alpha<1 $, $ N>2\alpha $, $ 1<q<\min\{2,2^*_\alpha-1\} $, $\lambda>0 $ is real parameter and $ 2^*_\alpha=\frac{2N}{N-2\alpha} $ is the fractional critical Sobolev exponent. Here $(-\Delta)^\alpha $ is the fractional Laplacian defined, up to a normalization constant,  as
\beq\label{1.2}
(-\Delta)^\alpha u(x)=P.V.\int_{\R^N}\frac{u(x)-u(y)}{|x-y|^{N+2\alpha}}dy,
\eqq
for $ x\in \R^N $, where P.V. denotes the principal value of the integral.

Concerning the weight functions $ a_\lambda(x) $ and $ b(x) $, we may assume that\\

$(H_1)$  $ a_\lambda=\lambda a_++a_-$, with $ a_{\pm}=\pm\max\{\pm a,0\}\not\equiv0 $, and $ b $ are continuous

\quad\quad in $ \bar{\Omega} $;\\

$(H_2)$  There exists a nonempty closed set
\[
M=\left\{x\in \bar{\Omega}\,\,\bigg\vert\,\,b(x)=\max_{\bar{\Omega}} b\equiv 1\right\}\subset \Omega
\]
\quad\quad \quad \quad and  a positive number $ t>\frac{N-2\alpha}{2} $ such that
\[
b(z)-b(x)=O(|x-z|^t)\quad
\]
\quad\quad \quad \quad holds uniformly for $z\in M$ in  the limit $x\rightarrow z$.\\

\begin{remark}\label{r1}
Let 
\[M_r=\{x\in \R^N\,|\,{\rm dist}(x,M)<r\}\quad for \,\, r>0.\]
By $(H_2)$, we may then assume that there exist constants $\eta_0, D_0$ and $r_0$ such that
\[
b(x)\ge \eta_0\quad for\,\,all\,\,x\in M_{r_0}\subset \Omega.
\]
and
\[
b(z)-b(x)\le D_0|x-z|^{t}\quad for\,\,all\,\,x\in B_{r_0}(z) \,\,and\,\,for\,\,all\,\, z\in M.
\]

\end{remark}

When $a_\lambda\equiv \lambda$ and $b\equiv1$, problem (\ref{1.1}) has been studied by Barrios et al. in \cite{BCSS}. They proved that there exists a positive $\Lambda$ such that
(\ref{1.1}) admits at least two solutions if $\lambda\in (0,\Lambda)$. One can also define a fractional power of the Laplacian using spectral decomposition. Problem (\ref{1.1}) for the spectral
factional Laplacian  has been treated in \cite{BCPS}. In this article, we study problem (\ref{1.1}) with sign-changing weight functions. Our first main result is

\begin{teo}\label{t1.1}
Suppose that $(H_1)$ and $(H_2)$ hold. Let
\[
\Lambda_0=\frac{q}{2}\cdot \frac{S_\alpha^{\frac{N(2-q)}{4\alpha}+\frac{q}{2}}}{\|a_+\|_{L^{q^*}(\Omega)}}\cdot \left(\frac{2-q}{2^*_\alpha-q}\right)^{\frac{(N-2\alpha)(2-q)}{4\alpha}}\cdot \left(\frac{2^*_\alpha-2}{2^*_\alpha-q}\right),
\]
where $ S_\alpha $ is the best Sobolev constant for the embedding of $ H^\alpha(\R^N) $ into $ L^{2^*_\alpha}(\R^N) $ (see (\ref{2.1}) below) and $ q^*=2^*_\alpha/(2^*_\alpha-q) $. Then
problem $(1.1)$ has at least two positive solutions if  $\lambda\in(0,\Lambda_0)$.
\end{teo}

We use variational methods to find positive solutions of equation (\ref{1.1}).
We denote by $H^\alpha(\R^N)$ the usual fractional Sobolev space endowed with the so-called Gagliardo norm
\beq\label{1.3}
\|u\|_{H^\alpha(\R^N)}=\|u\|_{L^2(\R^N)}+\left(\int_{\R^{2N}}\frac{|u(x)-u(y)|^2}{|x-y|^{N+2\alpha}}dxdy\right)^{1/2}.
\eqq
while $X_0^\alpha(\Omega)$ is the function space defined as
\beq\label{1.4}
X_0^\alpha(\Omega)=\{u\in H^\alpha(\R^N): u=0\,\,{\rm a.e.}\,\,{\rm in}\,\,\R^N\setminus\Omega\}.
\eqq
In $X_0^\alpha(\Omega)$ we consider the following norm
\beq\label{1.5}
\|u\|_{X_0^\alpha(\Omega)}=\left(\int_{\R^{2N}}\frac{|u(x)-u(y)|^2}{|x-y|^{N+2\alpha}}dxdy\right)^{1/2}.
\eqq
We also recall that $ (X_0^\alpha(\Omega), \|\cdot\|_{X_0^\alpha(\Omega)}) $ is a Hilbert space with scalar product
\beq\label{1.6}
\langle u,v\rangle_{X_0^\alpha(\Omega)}=\int_{\R^{2N}}\frac{(u(x)-u(y))(v(x)-v(y))}{|x-y|^{N+2\alpha}}dxdy,
\eqq
see Lemma 7 in \cite{SV}.
Note that by Proposition 3.6 in \cite{DPV} we have the following identities, up to constants,
\[
\|u\|_{X_0^\alpha(\Omega)}=\||\xi|^\alpha\mathcal{F}u\|_{L^2(\R^N)}=\|\mathcal{F}(-\Delta)^{\alpha/2} u\|_{L^2(\R^N)}=\left(\int_{\R^N}|(-\Delta)^{\alpha/2} u|^2dx\right)^{1/2}.
\]
We have used that if $u$ and $v$ in $X_0^\alpha(\Om)$, then
\[
\int_\Om v(-\Delta)^\alpha u dx =\int_{\R^N} (-\Delta)^{\alpha/2} v(-\Delta)^{\alpha/2} udx
\]
which yields the following definition

\begin{definition}\label{d1}
We say that $u\in X_0^\alpha(\Omega)$ is a weak solution of (\ref{1.1}) if for every $\varphi\in X_0^\alpha(\Omega)$,
one has
\beq\label{1.7}
\int_{\R^{2N}}\frac{(u(x)-u(y))(\varphi(x)-\varphi(y))}{|x-y|^{N+2\alpha}}dxdy=\int_\Omega a_\lambda|u|^{q-1}u\varphi dx+\int_\Omega b|u|^{p-1}u\varphi dx.
\eqq
\end{definition}

In this sequel we will omit the term weak when referring to solutions that satisfy the conditions of Definition \ref{d1}.
Associated with equation (\ref{1.1}), we consider the energy  functional $\Phi_\lambda$ in $X_0^\alpha(\Omega)$,
\[
\Phi_\lambda(u)=\frac{1}{2}\|u\|^2_{X_0^\alpha(\Omega)}-\frac{1}{q}\int_\Omega a_\lambda|u|^{q} dx-\frac{1}{2^*_\alpha}\int_\Omega b|u|^{2^*_\alpha} dx.
\]
As it is well known, when one uses the variational methods to find the critical points of the functional, some geometry structures are needed such as the mountain pass structure, the linking structures and so on. For problem (\ref{1.1}), the main difficulty lies in the functional may not posses such structures since the sign-changing weight. In order to overcome this difficulty, we turn to another approach, that is, the Nehari manifold, which was introduced by Nehari in \cite{N} and has been widely used in the literature, for example \cite{T,B,W1,W2,W3} and references therein for Laplace operator and also \cite{CDS,Y} for  the fractional Laplacian. The main idea of these articles lies in dividing the Nehari manifold into three parts and considering the infima of the functional on each part. More precise, the Nehari manifold for $ \Phi_\lambda(u) $ is defined as
\begin{eqnarray*}
\mathcal{N}_\lambda &=&\{u\in X_0^\alpha(\Omega)\,\,|\,\,\langle \Phi_\lambda^\prime(u),u\rangle=0\}\\
                    &=&\left\{u\in X_0^\alpha(\Omega)\,\,\bigg\vert\,\,\|u\|^2_{X_0^\alpha(\Omega)}-\int_\Omega a_\lambda|u|^{q} dx-\int_\Omega b|u|^{2^*_\alpha} dx=0\right\}.
\end{eqnarray*}
It is clear that all critical points of $ \Phi $ must be lie on $\mathcal{N}_\lambda  $, as we will see below, local minimizers on $ \mathcal{N}_\lambda $ are usually critical points of $ \Phi_\lambda $. By consider the fibering map $ h_u(t)= \Phi_\lambda(tu)$, we can divide that $\mathcal{N}_\lambda  $ into three subsets $ \mathcal{N}_\lambda^+, \mathcal{N}_\lambda^- $ and $ \mathcal{N}_\lambda^0 $ which correspond to local minima, local maxima and points of inflexion of fibbering maps. Then we can find that $\mathcal{N}_\lambda^0=\emptyset$ if $\lambda\in (0,\Lambda_0)$ and meanwhile there exists at least one positive solution in $ \mathcal{N}_\lambda^+$ and $ \mathcal{N}_\lambda^- $ respectively. Moreover,  by applying the Ljusternik-Schnirelmann category (see for example \cite{J}), we can show another multiplicity result. We would like point out that, if Y is a closed subset of a topological space $X$, the Lusternik-Schnirelman category $cat_X (Y )$ is the least number of closed and contractible sets in $X$ which cover $Y$ .  Here and in what follows, we denote $cat$ as the Ljusternik-Schnirelmann category. Recalling the definition of $M$ and $M_\delta$ in $(H_2)$ and Remark \ref{r1} respectively and using  the Ljusternik-Schnirelmann category, we can prove that

\begin{teo}\label{t1.2}
Suppose that $(H_1)$ and $(H_2)$ hold. For each $ \delta<r_0 $ (see Remark \ref{r1}), then there exists $ 0<\Lambda_\delta \le \Lambda_0$ such that problem (\ref{1.1}) has at least $ cat_{M_\delta}(M)+1 $ positive solutions for each $ \lambda\in (0,\Lambda_\delta) $.
\end{teo}

When $\alpha=1$ and $u=0$ on $\partial \Omega$, de Pavia \cite{P} studied  sufficient small $\lambda$ and obtained a globalized result, indicating that there exists a $\lambda^*$ such that (\ref{1.1}) has at least two solutions if $\lambda\in (0,\lambda^*)$.
In \cite{P}, they requires that one of the weight functions is non-negative with a non-empty domain for which $a(x)$ and $b(x)$ are both positive. 
In order to overcome the nonnegative assumptions on the weight functions,  Chen et al. \cite{CW} recently by studying the decomposition of the Nehari manifold relaxed the conditions  of the weight functions set out by de Pavia \cite{P} with hypotheses $(H_1)-(H_2)$ (without imposing the non-negativity constraint on the weight functions  $a(x)$ and $b(x)$) and investigate the solution structure of (\ref{1.1}). This method is also used in \cite{W1,W2,W3,T,B} and reference therein. Furthermore, in \cite{CW} the authors also proved there exists at least $cat_{M_\delta}(M)+1 $ positive solutions based on the concentration-compactness principle and the Lusternik-Schnirelman category.  The concentration-compactness principle for the fractional Laplacian
is obtained by Palatucci and Pisante \cite{PP} recently. Thus, we would like to extend the result in \cite{CW} to equation (\ref{1.1}), Theorem \ref{t1.2}.

This article is organized as follows. In Section 2 we give some notations and  preliminaries for the Nehari manifold. Sections 3 and 4 are devoted to prove the multiplicity of positive
 solutions of equation (\ref{1.1}), Theorem \ref{t1.1} and Theorem \ref{t1.2}, respectively.

\setcounter{equation}{0}
\section{Preliminaries}

We start this section by recalling the best  Sobolev constant $ S_\alpha $ for the embedding of $ H^\alpha(\R^N) $ into $ L^{2^*_\alpha}(\R^N) $, which is defined as
\beq\label{2.1}
S_\alpha=\inf_{H^\alpha(\R^N)\setminus\{0\}}\frac{\int_{\R^{2N}}\frac{|u(x)-u(y)|^2}{|x-y|^{N+2\alpha}}dxdy}{\left(\int_{\R^N}|u|^{2^*_\alpha}dx\right)^{2/2^*_\alpha}}>0.
\eqq
By Theorem 1.1 in \cite{CT}, the infimun in (\ref{2.1}) is attained at the function
\beq\label{2.2}
u_0(x)=\kappa/(|x-x_0|^2+\mu^2)^{(N-2\alpha)/2}
\eqq
where $ \kappa\in\R $, $ \mu>0 $ and $ x_0\in\R^N $ are fixed constants. Moreover,  let
\[
\tilde{u}(x)=u_0(x/S_\alpha^{1/2\alpha})\quad{\rm for}\,\,x\in\R^N,
\]
then $ \tilde{u} $ is a positive solution of the critical problem
\beq\label{2.3}
(-\Delta)^\alpha u=|u|^{2^*_\alpha-1}u\quad{\rm in}\,\,\R^N.
\eqq
Furthermore, for any $ \varepsilon>0 $, we define
\beq\label{2.17}
U_\varepsilon(x)=\varepsilon^{-\frac{N-2\alpha}{2}}\tilde{u}\left(x/\varepsilon\right)\quad{\rm for}\,\,x\in\R^N,
\eqq
then $ U_\varepsilon $ satisfying (\ref{2.3}) and also
\[
\|U_\varepsilon\|_{H^\alpha(\R^N)}^2=\|U_\varepsilon\|_{L^{2^*_\alpha}(\R^N)}^{2^*_\alpha}=S_\alpha^{N/2\alpha}.
\]

We define the Palais-Smale (PS)-sequences and (PS)-condition in $X_0^\alpha(\Omega)$ for $\Phi_\lambda$ as follows.

\begin{definition}\label{d2.1}
(1) For $c\in\R$, a sequence $\{u_n\}$ is a $(PS)_c$-sequence in $X_0^\alpha(\Omega)$ for $\Phi_\lambda$ if $\Phi_\lambda(u_n)=c+o(1)$ and $\Phi_\lambda^\prime(u_n)=o(1)$ strongly in $(X_0^\alpha(\Omega))^*$ as $n\rightarrow\infty$.

(2) $\Phi_\lambda$ satisfies the $(PS)_c$-condition in $X_0^\alpha(\Omega)$ if every $(PS)_c$-sequence in $X_0^\alpha(\Omega)$ for  $\Phi_\lambda$
contains a convergent subsequence.
\end{definition}

Since the energy functional $\Phi_\lambda$ is not bounded below on $X_0^\alpha(\Omega)$, it is useful to consider the functional on the Nehari manifold
$\mathcal{N}_\lambda$.
Moreover, we have the following result.

\begin{lemma}\label{l2.1}
The energy functional $\Phi_\lambda$ is coercive and bounded below on $\mathcal{N}_\lambda$.
\end{lemma}

{\bf Proof.} By H\"older and Sobolev inequalities, for $u\in \mathcal{N}_\lambda$, we have
\begin{eqnarray}\label{2.16}
\Phi_\lambda(u)&=&\left(\frac{1}{2}-\frac{1}{2^*_\alpha}\right)\|u\|^2_{X_0^\alpha(\Omega)}-\left(\frac{1}{q}-\frac{1}{2^*_\alpha}\right)\int_\Omega (\lambda a_++a_-)|u|^{q} dx\\\nonumber
&\ge&\left(\frac{1}{2}-\frac{1}{2^*_\alpha}\right)\|u\|^2_{X_0^\alpha(\Omega)}-\left(\frac{1}{q}-\frac{1}{2^*_\alpha}\right)\int_\Omega \lambda a_+|u|^{q} dx\\\label{2.15}
&\ge&\left(\frac{1}{2}-\frac{1}{2^*_\alpha}\right)\|u\|^2_{X_0^\alpha(\Omega)}-\lambda\left(\frac{1}{q}-\frac{1}{2^*_\alpha}\right)\|a_+\|_{L^{q^*}(\Omega)}S_\alpha^{-\frac{q}{2}}\|u\|_{X_0^\alpha(\Omega)}^q,
\end{eqnarray}
where $q^*=2^*_\alpha/(2^*_\alpha-q)$.
Then $\Phi_\lambda$  is coercive and bounded below on $\mathcal{N}_\lambda$. $\Box$\\

The Nehari manifold $ \mathcal{N}_\lambda $ is closely related to the behaviour of the function of the form $ h_u :t\rightarrow\Phi_\lambda(tu)$ for $ t>0 $. Such map are know as fibering maps that dates back to the fundamental works \cite{P1,C,R,DP}. If $ u\in X_0^\alpha(\Omega) $, we have
\begin{eqnarray*}
h_u(t)&=&\frac{t^2}{2}\|u\|^2_{X_0^\alpha(\Omega)}-\frac{t^q}{q}\int_\Omega a_\lambda|u|^{q} dx-\frac{t^{2^*_\alpha}}{2^*_\alpha}\int_\Omega b|u|^{2^*_\alpha} dx;\\
h^\prime _u(t)&=&t\|u\|^2_{X_0^\alpha(\Omega)}-t^{q-1}\int_\Omega a_\lambda|u|^{q} dx-t^{2^*_\alpha-1}\int_\Omega b|u|^{2^*_\alpha} dx;\\
h^{\prime\prime} _u(t)&=&\|u\|^2_{X_0^\alpha(\Omega)}-(q-1)t^{q-2}\int_\Omega a_\lambda|u|^{q} dx-(2^*_\alpha-1)t^{2^*_\alpha-2}\int_\Omega b|u|^{2^*_\alpha} dx.
\end{eqnarray*}
We observe that
\[
h^\prime _u(t)=\langle\Phi_\lambda^\prime(tu),u\rangle=\frac{1}{t}\langle\Phi_\lambda^\prime(tu),tu\rangle
\]
and thus, for $ u\in  X_0^\alpha(\Omega)\setminus\{0\} $ and $ t>0 $, $ h^\prime _u(t)=0 $ if and only if $ tu\in  \mathcal{N}_\lambda  $, that is, positive critical points of $ h_u $ correspond points on the Nehari manifold. In particular, $ h^\prime _u(1)=0 $ if and only if
$ u\in  \mathcal{N}_\lambda  $. So it is natural to split $\mathcal{N}_\lambda  $ into three parts corresponding local minimal, local maximum and points of inflection. Accordingly, we define
\begin{eqnarray*}
\mathcal{N}_\lambda^+&=&\{u\in \mathcal{N}_\lambda\,\,|\,\,h^{\prime\prime} _u(1)>0\};\\
\mathcal{N}_\lambda^0&=&\{u\in \mathcal{N}_\lambda\,\,|\,\,h^{\prime\prime} _u(1)=0\};\\
\mathcal{N}_\lambda^-&=&\{u\in \mathcal{N}_\lambda\,\,|\,\,h^{\prime\prime} _u(1)<0\}.
\end{eqnarray*}

Next, we establish some basic properties of $ \mathcal{N}_\lambda^+, \mathcal{N}_\lambda^0 $, and $ \mathcal{N}_\lambda^- $.

\begin{lemma}\label{l2.2}
Suppose that $u_0$ is a local minimizer of $\Phi_\lambda$ on $\mathcal{N}_\lambda$ and $u_0\not\in \mathcal{N}_\lambda^0$.
Then $\Phi_\lambda^\prime(u_0)=0$  in $(X_0^\alpha(\Omega))^*$, where  $(X_0^\alpha(\Omega))^*$ is the dual space of  $X_0^\alpha(\Omega)$.
\end{lemma}

{\bf Proof.} If $ u_0 $ is a local minimizer for $\Phi_\lambda$ on $\mathcal{N}_\lambda$, then $ u_0 $ is a solution of the optimization problem
\[
{\rm minimizer}\,\, \Phi_\lambda(u) \,\,{\rm subject \,\,to} \,\,J(u)=0,
\]
where $ J(u)= \|u\|^2_{X_0^\alpha(\Omega)}-\int_\Omega a_\lambda|u|^{q} dx-\int_\Omega b|u|^{2^*_\alpha} dx$. Hence, by the theory of Lagrange multipliers, there exists $ \mu\in\R $  such that $ \Phi_\lambda^\prime(u_0)=\mu J^\prime(u_0) $. Thus we have
\beq\label{2.18}
\langle\Phi_\lambda^\prime(u_0), u_0\rangle=\mu \langle J^\prime(u_0), u_0\rangle.
\eqq
Since $ u_0\in \mathcal{N}_\lambda $, we have that $ \|u_0\|^2_{X_0^\alpha(\Omega)}-\int_\Omega a_\lambda|u_0|^{q} dx-\int_\Omega b|u_0|^{2^*_\alpha} dx=0$. Hence,
\begin{eqnarray*}
\langle J^\prime(u_0), u_0\rangle&=&2\|u_0\|^2_{X_0^\alpha(\Omega)}-q\int_\Omega a_\lambda|u_0|^{q} dx-2^*_\alpha\int_\Omega b|u_0|^{2^*_\alpha} dx\\
&=&\|u_0\|^2_{X_0^\alpha(\Omega)}-(q-1)\int_\Omega a_\lambda|u_0|^{q} dx-(2^*_\alpha-1)\int_\Omega b|u_0|^{2^*_\alpha} dx.
\end{eqnarray*}
So, if $u_0\not\in \mathcal{N}_\lambda^0$, $\langle J^\prime(u_0), u_0\rangle\not=0$ and thus $ \mu=0 $ by (\ref{2.18}). Hence, we complete the proof. $\Box$\\

For each $ u\in \mathcal{N}_\lambda $, we know that
\begin{eqnarray}\nonumber
h^{\prime\prime} _u(1)&=&\|u\|^2_{X_0^\alpha(\Omega)}-(q-1)\int_\Omega a_\lambda|u|^{q} dx-(2^*_\alpha-1)\int_\Omega b|u|^{2^*_\alpha} dx\\\label{2.7}
&=&(2-2^*_\alpha)\|u\|^2_{X_0^\alpha(\Omega)}-(q-2^*_\alpha)\int_\Omega a_\lambda|u|^{q} dx\\\label{2.8}
&=&(2-q)\|u\|^2_{X_0^\alpha(\Omega)}-(2^*_\alpha-q)\int_\Omega b|u|^{2^*_\alpha} dx.
\end{eqnarray}
Then we have following result.

\begin{lemma}\label{l2.3}
(1) For any $ u\in \mathcal{N}_\lambda^+\cup \mathcal{N}_\lambda^0 $, we have $ \int_\Omega a_\lambda|u|^{q} dx>0 $;\\
(2) For any $ u\in \mathcal{N}_\lambda^-$, we have $ \int_\Omega b|u|^{2^*_\alpha} dx>0 $.
\end{lemma}

{\bf Proof.}  By the definitions of $\mathcal{N}_\lambda^+$ and $\mathcal{N}_\lambda^0$, it is easy to get that $ \int_\Omega a_\lambda|u|^{q} dx>0 $  from (\ref{2.7}). Similarly, the definition of $\mathcal{N}_\lambda^-$ and (\ref{2.8}) imply that $ \int_\Omega b|u|^{2^*_\alpha} dx>0 $. $\Box$\\

Let $ \Lambda_1= \frac{S_\alpha^{\frac{N(2-q)}{4\alpha}+\frac{q}{2}}}{\|a_+\|_{L^{q^*}(\Omega)}}\cdot \left(\frac{2-q}{2^*_\alpha-q}\right)^{\frac{(N-2\alpha)(2-q)}{4\alpha}}\cdot \left(\frac{2^*_\alpha-2}{2^*_\alpha-q}\right)$. Then we have the following result.

\begin{lemma}\label{l2.4}
We have $ \mathcal{N}_\lambda^0=\emptyset $ for all $ \lambda<\Lambda_1 $.
\end{lemma}

{\bf Proof.} We prove it by contradiction arguments. Suppose that there exists  $ \lambda<\Lambda_1 $ such that $ \mathcal{N}_\lambda^0\not=\emptyset $. Then, for $ u_0\in \mathcal{N}_\lambda^0 $, by (\ref{2.7}) and the H\"older and Sobolev inequalities, we have
\begin{eqnarray*}
\|u\|^2_{X_0^\alpha(\Omega)}&=&\frac{2^*_\alpha-q}{2^*_\alpha-2}\int_\Omega a_\lambda|u|^{q} dx\\
&\le &\frac{2^*_\alpha-q}{2^*_\alpha-2}\int_\Omega \lambda a_+|u|^{q}\\
&\le &\lambda\cdot\frac{2^*_\alpha-q}{2^*_\alpha-2}\|a_+\|_{L^{q^*}(\Omega)}\|u\|^{q}_{L^{2^*_\alpha}(\Omega)}\\
&\le &\lambda\cdot\frac{2^*_\alpha-q}{2^*_\alpha-2}\|a_+\|_{L^{q^*}(\Omega)}S_\alpha^{-\frac{q}{2}}\|u\|^{q}_{X_0^\alpha(\Omega)}
\end{eqnarray*}
and so
\beq\label{2.9}
\|u\|^{2-q}_{X_0^\alpha(\Omega)}\le \lambda\cdot\frac{2^*_\alpha-q}{2^*_\alpha-2}\|a_+\|_{L^{q^*}(\Omega)}S_\alpha^{-\frac{q}{2}}.
\eqq
Similarly, by (\ref{2.8}) the H\"older and Sobolev inequalities, we have
\beq\label{2.10}
\|u\|_{X_0^\alpha(\Omega)}\ge \left(\frac{2-q}{2^*_\alpha-q}\right)^{\frac{(N-2\alpha)}{4\alpha}}S_\alpha^{\frac{N}{4\alpha}}
\eqq
since $ \max_{\bar{\Omega}}b(x)\equiv 1 $.

Hence, combining (\ref{2.9}) and (\ref{2.10}), we must have
\[
\lambda\ge \frac{S_\alpha^{\frac{N(2-q)}{4\alpha}+\frac{q}{2}}}{\|a_+\|_{L^{q^*}(\Omega)}}\cdot \left(\frac{2-q}{2^*_\alpha-q}\right)^{\frac{(N-2\alpha)(2-q)}{4\alpha}}\cdot \left(\frac{2^*_\alpha-2}{2^*_\alpha-q}\right)=\Lambda_1,
\]
which is a contradiction. This completes the proof. $ \Box $\\

In order to get a better understanding of the Nehari manifold and the fibering maps, we considering the function $ m_u:\R^+\rightarrow \R $ defined by
\beq\label{2.11}
m_u(t)=t^{2-q}\|u\|_{X_0^\alpha(\Omega)}^2-t^{2^*_\alpha-q}\int_\Omega b|u|^{2^*_\alpha} dx\quad {\rm for}\,\,t>0.
\eqq
It is clear that $ tu\in \mathcal{N}_\lambda $ if and only $ m_u(t)= \int_\Omega a_\lambda|u|^{q}$. Moreover,
\beq\label{2.12}
m_u^\prime(t)=(2-q)t^{1-q}\|u\|_{X^0_\alpha(\Omega)}^2-(2^*_\alpha-q)t^{2^*_\alpha-q-1}\int_\Omega b|u|^{2^*_\alpha} dx
\eqq
and it is easy to see that, if $ tu\in \mathcal{N}_\lambda $, then $ t^{q-1}m^\prime_u(t)=h_u^{\prime\prime}(t) $. Hence $ tu\in  \mathcal{N}_\lambda ^+ $ (or $  \mathcal{N}_\lambda ^- $) if and only if $ m^\prime_u(t)>0 $ (or $ <0 $).

For every $u\in X_0^\alpha(\Omega)\setminus\{0\}$ with $\int_\Omega b|u|^{2^*_\alpha} dx>0$, we let
\beq\label{2.6}
t_{{ \rm max}}(u)=\left(\frac{(2-q)\|u\|_{X_0^\alpha(\Omega)}^2}{(2^*_\alpha-q)\int_\Omega b|u|^{2^*_\alpha} dx}\right)^{\frac{N-2\alpha}{4\alpha}}>0,
\eqq
which leads the following lemma.

\begin{lemma}\label{l2.5}
Suppose that $\int_\Omega b|u|^{2^*_\alpha} dx>0$. Then for each $u\in X_0^\alpha(\Omega)\setminus\{0\}$ and $\lambda\in (0,\Lambda_1)$, we have that

(1) if $\int_\Omega a_\lambda|u|^{q} dx\le0$, then there exists a unique $t^-=t^-(u)>t_{{\rm max}}(u)$ such that
$t^-u\in \mathcal{N}_\lambda^-$ and
\beq\label{2.13}
\Phi_\lambda(t^-u)=\sup_{t\ge0}\Phi_\lambda(tu).
\eqq

(2) if $\int_\Omega a_\lambda|u|^{q} dx>0$, then there exists a unique $0<t^+=t^+(u)<t_{{\rm max}}(u)<t^-$ such that
$t^+u\in \mathcal{N}_\lambda^+$,  $t^-u\in \mathcal{N}_\lambda^-$ and
\beq\label{2.14}
\Phi_\lambda(t^+u)=\inf_{0\le t\le t_{{\rm max}}(u)}\Phi_\lambda(tu),\quad \Phi_\lambda(t^-u)=\sup_{t\ge t^+}\Phi_\lambda(tu).
\eqq
\end{lemma}

{\bf Proof.} By (\ref{2.12}), we know $ t_{{\rm max}} $ is the unique critical point of $ m_u $ and $ m_u $ is strictly increasing on $ (0, t_{{\rm max}}) $  and strictly decreasing on $ ( t_{{\rm max}},\infty) $ with $ \lim_{t\rightarrow\infty}m_u(t)=-\infty $. Moreover, by the H\"older and Sobolev inequalities, we have that
\begin{eqnarray*}
m_u(t_{{\rm max}})
&\ge&\left(\frac{2-q}{2^*_\alpha-q}\right)^{\frac{(N-2\alpha)(2-q)}{4\alpha}}\cdot \left(\frac{2^*_\alpha-2}{2^*_\alpha-q}\right) S_\alpha^{\frac{N(2-q)}{4\alpha}+\frac{q}{2}}\|u\|_{L^{2^*_\alpha}(\Omega)}^q\\
&\ge&\left(\frac{2-q}{2^*_\alpha-q}\right)^{\frac{(N-2\alpha)(2-q)}{4\alpha}}\cdot \left(\frac{2^*_\alpha-2}{2^*_\alpha-q}\right) S_\alpha^{\frac{N(2-q)}{4\alpha}+\frac{q}{2}}\frac{\int_\Omega a_\lambda|u|^q dx}{\lambda\|a_+\|_{L^{q^*}(\Omega)}}\\
&=& \Lambda_1\lambda^{-1}\int_\Omega a_\lambda|u|^q dx>\int_\Omega a_\lambda|u|^q dx.
\end{eqnarray*}

Next, we fix $u\in X_0^\alpha(\Omega)\setminus\{0\}$. Suppose that $\int_\Omega a_\lambda|u|^{q} dx\le0$. Then $ m_u(t)= \int_\Omega a_\lambda|u|^{q}$ has unique solution $t^->t_{{\rm max}}$ and $ m_u^\prime(t^-)<0 $. Hence $ h_u $ has a unique turning point at $ t=t^- $ and $ h^{\prime\prime}(t^-)<0 $. Thus $t^-u\in \mathcal{N}_\lambda^-$ and (\ref{2.13}) holds.

Suppose $\int_\Omega a_\lambda|u|^{q} dx>0$. Since $ m_u(t_{{\rm max}}) >\int_\Omega a_\lambda\|u\|^q dx$, the equation $ m_u(t)= \int_\Omega a_\lambda|u|^{q}$  has exactly two solutions
$0<t^+<t_{{\rm max}}(u)<t^-$ such that $ m_u^\prime(t^+)>0 $ and $ m_u^\prime(t^-)<0 $. Hence, there are two multiplies of $ u $ lying in $  \mathcal{N}_\lambda$, that is, $t^+u\in \mathcal{N}_\lambda^+$ and $t^-u\in \mathcal{N}_\lambda^-$. Thus $ h_u $ has turning points at $ t=t^+ $ and $ t=t^- $ with $ h^{\prime\prime}(t^+)<0 $ and $ h^{\prime\prime}(t^-)<0 $. Thus, $ h_u $ is decreasing on $ (0,t^+) $, increasing on $ (t^-,t^+) $ and decreasing on $ (t^-,\infty) $. Hence $ (\ref{2.14}) $ holds. $\Box$\\

\setcounter{equation}{0}
\section{Proof of Theorem \ref{t1.1}}

In this section, we prove Theorem \ref{t1.1} by variational methods. Firstly, by Lemma \ref{l2.5}, we know that $ \mathcal{N}_\lambda^+ $ and $ \mathcal{N}_\lambda^- $ are non-empty. Moreover, by Lemma \ref{l2.4}, we can write $ \mathcal{N}_\lambda=\mathcal{N}_\lambda^+\cup \mathcal{N}_\lambda^- $ and Lemma \ref{l2.1}, we can define
\[
c_\lambda^+=\inf_{u\in \mathcal{N}_\lambda^+}\Phi_\lambda(u)\quad{\rm and}\quad c_\lambda^-=\inf_{u\in \mathcal{N}_\lambda^-}\Phi_\lambda(u).
\]
\begin{lemma}\label{l3.1}

(1) For all $\lambda\in (0,\Lambda_1)$, we have $c_\lambda^+<0$;

(2) If $\lambda<\Lambda_0=\frac{1}{2}q\Lambda_1$, then $c_\lambda^->0$. In particular, $c_\lambda^+=\inf_{u\in \mathcal{N}_\lambda}\Phi_\lambda(u)$ for all $\lambda\in (0,\Lambda_0)$.
\end{lemma}

{\bf Proof.} (1) Let $ u\in \mathcal{N}_\lambda^+ $. Then, by (\ref{2.7}), we have
\[
\|u\|_{X_0^\alpha(\Omega)}^2<\frac{2^*_\alpha-q}{2^*_\alpha-2}\int_\Omega a_\lambda|u|^q dx.
\]
Hence, by (\ref{2.16}) and Lemma \ref{l2.3}, we have
\begin{eqnarray*}
\Phi_\lambda(u)&=&\left(\frac{1}{2}-\frac{1}{2^*_\alpha}\right)\|u\|^2_{X_0^\alpha(\Omega)}-\left(\frac{1}{q}-\frac{1}{2^*_\alpha}\right)\int_\Omega  a_\lambda|u|^{q} dx\\
&<&-\frac{(2^*_\alpha-q)(2-q)}{2q\cdot 2^*_\alpha}\int_\Omega  a_\lambda|u|^{q} dx<0.
\end{eqnarray*}
Thus, $c_\lambda^+<0$.

(2) Let $ u\in \mathcal{N}_\lambda^- $. Then, by (\ref{2.8}), we have
\[
\frac{2-q}{2^*_\alpha-q}\|u\|_{X_0^\alpha(\Omega)}^2<\int_\Omega b|u|^{2^*_\alpha} dx\le\int_\Omega |u|^{2^*_\alpha} dx\le S_\alpha^{-\frac{2^*_\alpha}{2}}\|u\|_{X_0^\alpha(\Omega)}^{2^*_\alpha}
\]
and so
\[
\|u\|_{X_0^\alpha(\Omega)}> S_\alpha^{\frac{N}{4\alpha}}\left(\frac{2-q}{2^*_\alpha-q}\right)^{\frac{N-2\alpha}{4\alpha}}.
\]
Therefore, by (\ref{2.15}), we know
\begin{eqnarray*}
\Phi_\lambda(u)
&\ge&\left(\frac{1}{2}-\frac{1}{2^*_\alpha}\right)\|u\|^2_{X_0^\alpha(\Omega)}-\lambda\left(\frac{1}{q}-\frac{1}{2^*_\alpha}\right)\|a_+\|_{L^{q^*}(\Omega)}S_\alpha^{-q/2}\|u\|_{X_0^\alpha(\Omega)}^q\\
&>& \|u\|^q_{X_0^\alpha(\Omega)}\left(\frac{\alpha}{N}S_\alpha^{\frac{N(2-q)}{4\alpha}}\left(\frac{2-q}{2^*_\alpha-q}\right)^{\frac{(N-2\alpha)(2-q)}{4\alpha}}-\lambda\frac{q-2^*_\alpha}{q2^*_\alpha}\|a_+\|_{L^{q^*}(\Omega)}S_\alpha^{-\frac{q}{2}}\right).
\end{eqnarray*}
Thus, if $ \lambda<\Lambda_0=\frac{q}{2}\Lambda_1$, then $ \Phi_\lambda(u)>0 $. This completes the proof. $\Box$\\

We need the following proposition for the precise description of the Palais-Smale sequence of $\Phi_\lambda$.

\begin{proposition}\label{p3.1}
Each sequence $ \{u_n\}\subset \mathcal{N}_\lambda $ that satisfies

(1) $ \Phi_\lambda(u_n)=c+o(1) $ with $ c<c_\lambda^++\frac{\alpha}{N}S_\alpha^{\frac{N}{2\alpha}} $;

(2) $ \Phi^\prime_\lambda(u_n)=o(1) $ in $ (X_0^\alpha(\Omega))^*$ \\
has a convergent subsequence.
\end{proposition}

The proof of Proposition \ref{p3.1} is very similar to  Proposition 3.2 in \cite{W3}, we omit it here.

Next, we establish the existence of a local minimum for $\Phi_\lambda$ on $\mathcal{N}_\lambda^+$.

\begin{teo}\label{t3.1}
For each $0<\lambda<\Lambda_0$, the functional $\Phi_\lambda$ has a minimizer $u_\lambda^+$ in $\mathcal{N}_\lambda^+$
satisfying that\\

(1) $\Phi_\lambda(u_\lambda^+)=c_\lambda^+=\inf_{u\in \mathcal{N}_\lambda^+}\Phi_\lambda(u)$;\\

(2) $u_\lambda^+$ is a positive solution of (\ref{1.1}).
\end{teo}

{\bf Proof.} By Lemma \ref{l2.1}, we know $ \Phi_\lambda $ is bounded blow on $ \mathcal{N}_\lambda $ as well as $  \mathcal{N}_\lambda^+ $. Thus, by Ekeland variational principle \cite{E}, there exists $ \{u_n\}\subset \mathcal{N}_\lambda^+ $ such that it is a $ (PS)_{c_\lambda^+} $-sequence for $ \Phi_\lambda $. Then by Proposition \ref{p3.1}, there exists a subsequence of $ \{u_n\} $ such that $u_n\rightarrow u_\lambda^+ $ strongly in $ X_0^\alpha(\Omega) $. Moreover,
$\Phi_\lambda(|u_\lambda^+|)\le \Phi_\lambda(u_\lambda^+)$ (see (A.11) in \cite{SV1}) and $|u_\lambda^+|\in \mathcal{N}_\lambda^+$, by Lemma \ref{l2.2}, we may assume $ u_\lambda^+ $ is a positive solution of (\ref{1.1}). $ \Box $ \\

Next, we consider a cut-off function $ \eta\in C^\infty(\R^N) $ with $ 0\le\eta\le1 $, $ |\nabla\eta|\le C $, $ \eta=1 $ if $ |x|\le r_0/2 $ and  $ \eta=0 $ if $ |x|\ge r_0 $. For any $ z\in M $ (see hypothesis $(H_2)$), let
\beq\label{3.2}
w_{\varepsilon,z}(x)=\eta(x-z)U_\varepsilon(x-z),
\eqq
where $ U_\varepsilon $ given by (\ref{2.17}) with $ x_0=0 $. By similar argument as Propositions 21 and 22 in \cite{SV2}, we have that
\beq\label{3.3}
\|w_{\varepsilon,z}\|_{X_0^\alpha(\Omega)}^2=S_\alpha^{\frac{N}{2\alpha}}+O(\varepsilon^{N-2\alpha})\quad{\rm and}\quad \int_\Omega |w_{\varepsilon,z}|^{2_*^\alpha}dx=S_\alpha^{\frac{N}{2\alpha}}+O(\varepsilon^{N})
\eqq
hold uniformly for $ z\in M $. Hence, by using (\ref{3.3}) and taking a similar argument as Lemmas 3.1 and 3.2 in \cite{CW}, we can get the following estimates.

\begin{lemma}\label{l3.2}
(1) $\int_\Omega b |w_{\varepsilon,z}|^{2_*^\alpha}dx=S_\alpha^{\frac{N}{2\alpha}}+o(\varepsilon^{\frac{N-2\alpha}{2}})$ uniformly for $ z\in M $;

(2)$\int_\Omega |w_{\varepsilon,z}|^{q}dx=o(\varepsilon^{\frac{N-2\alpha}{2}})$ uniformly for $ z\in M $.
\end{lemma}

{\bf Proof.} (1) $w_{\varepsilon,z}$ is given by (\ref{3.2}). We define function $\tilde{b}:\R^N\rightarrow\R$ is an extension of $b$ by  $\tilde{b}(x)=b(x)$ if $x\in \bar{\Om}$ and $\tilde{b}(x)=0$ if $x\in \R^N\setminus\bar{\Om}$. 

By the definition of $U_\varepsilon$ (see (\ref{2.17})), we have
\begin{eqnarray*}
\int_{\Om}b|w_{\varepsilon,z}|^{2^*_\alpha}dx&=&\int_{B_{r_0}(z)}b(x)|\eta(x-z)U_\varepsilon(x-z)|^{2^*_\alpha}dx\\
&=&\int_{B_{r_0}(0)}b(x+z)|\eta(x)U_\varepsilon(x)|^{2^*_\alpha}dx\\
&=&\int_{B_{r_0}(0)}\frac{\ve^NS_\alpha^{\frac{N}{\alpha}}\kappa^{2^*_\alpha}}{\left(|x|^2+\ve^2S_\alpha^{\frac{1}{\alpha}}\mu^2\right)^N}\tilde{b}(x+z)\eta^{2^*_\alpha}(x)dx\\
&:=&\int_{B_{r_0}(0)}\frac{C_0\ve^N}{\left(|x|^2+C_1\ve^2\right)^N}\tilde{b}(x+z)\eta^{2^*_\alpha}(x)dx,\\
\end{eqnarray*}
where $C_0=S_\alpha^{\frac{N}{\alpha}}\kappa^{2^*_\alpha}$ and $C_1=S_\alpha^{\frac{1}{\alpha}}\mu^2$.

Next, by assumption $(H_2)$ and $b(z)=1$ since $z\in M$, we can see that
\begin{eqnarray*}
0&\le& \int_{\R^N}|U_\varepsilon|^{2^*_\alpha}dx-\int_{\Om}b|w_{\varepsilon,z}|^{2^*_\alpha}dx\\
&=&\int_{\R^N}\frac{C_0\ve^N}{\left(|x|^2+C_1\ve^2\right)^N}\left(1-\tilde{b}(x+z)\eta^{2^*_\alpha}(x)\right)dx\\
&=&\int_{\R^N\setminus B_{\frac{r_0}{2}}(0)}\frac{C_0\ve^N}{\left(|x|^2+C_1\ve^2\right)^N}\left(1-\tilde{b}(x+z)\eta^{2^*_\alpha}(x)\right)dx\\
&+&\int_{ B_{\frac{r_0}{2}}(0)}\frac{C_0\ve^N}{\left(|x|^2+C_1\ve^2\right)^N}\left(1-\tilde{b}(x+z)\eta^{2^*_\alpha}(x)\right)dx\\
&\le& C_0\ve^N \int_{\R^N\setminus B_{\frac{r_0}{2}}(0)}\frac{1}{|x|^{2N}}dx+D_0C_0\ve^N\int_{ B_{\frac{r_0}{2}}(0)}\frac{|x|^\rho}{\left(|x|^2+C_1\ve^2\right)^N}dx\\
&\le& C_0\ve^N \int^{\infty}_{\frac{r_0}{2}}r^{-N-1}dr+D_0C_0\ve^N\int_0^{\frac{r_0}{2}}\frac{r^{\rho+N-1}}{\left(r^2+C_1\ve^2\right)^N}dr\\
&\le& O(\ve^N)+D_0C_0\ve^N\int_0^{\ve}\frac{r^{\rho+N-1}}{\left(r^2+C_1\ve^2\right)^N}dr+D_0C_0\ve^N\int_{\ve}^{\frac{r_0}{2}}\frac{r^{\rho+N-1}}{\left(r^2+C_1\ve^2\right)^N}dr\\
&\le& O(\ve^N)+D_1C_0\ve^{-N}\int_0^{\ve}r^{\rho+N-1}dr+D_2C_0\ve^N\int_{\ve}^{\frac{r_0}{2}}r^{\rho-N-1}dr\\
&=&\left\{\begin{array}{l@{\quad }l}
O(\ve^N)+O(\ve^\rho) &{\rm if}\,\,\rho\not=N,\\
O(\ve^N)+C_2\ve^N|{\rm ln} \ve|\,\,&{\rm if}\,\,\rho=N.
 \end{array}
 \right.
\end{eqnarray*}
This implies that
\[
\int_\Omega b |w_{\varepsilon,z}|^{2_*^\alpha}dx=S_\alpha^{\frac{N}{2\alpha}}+o(\varepsilon^{\frac{N-2\alpha}{2}})
\]
uniformly for $z\in M$ since $\rho>(N-2\alpha)/2$ and $ \int_{\R^N}|U_\varepsilon|^{2^*_\alpha}dx=S_\alpha^{\frac{N}{2\alpha}}$.

(2) Since
\begin{eqnarray*}
\int_{\Om}|w_{\ve,z}|^qdx&=&\int_{B_{r_0}(0)}\eta^q(x)U_\ve^q(x)dx\\
&=&\int_{B_{r_0}(0)}\eta^q(x)\frac{C_0\ve^{\frac{q(N-2\alpha)}{2}}}{\left(|x|^2+C_1\ve^2\right)^{\frac{q(N-2\alpha)}{2}}}dx\\
&=&\int_{B_{\ve}(0)}\eta^q(x)\frac{C_0\ve^{\frac{q(N-2\alpha)}{2}}}{\left(|x|^2+C_1\ve^2\right)^{\frac{q(N-2\alpha)}{2}}}dx\\
&+&\int_{B_{r_0}(0)\setminus B_\ve(0)}\eta^q(x)\frac{C_0\ve^{\frac{q(N-2\alpha)}{2}}}{\left(|x|^2+C_1\ve^2\right)^{\frac{q(N-2\alpha)}{2}}}dx\\
&\le&C_3 \int_{B_{\ve}(0)}\frac{\ve^{\frac{q(N-2\alpha)}{2}}}{\ve^{q(N-2\alpha)}}dx+\int_{B_{r_0}(0)\setminus B_\ve(0)}\frac{C_0\ve^{\frac{q(N-2\alpha)}{2}}}{|x|^{q(N-2\alpha)}}dx\\
&=&\left\{\begin{array}{l@{\quad }l}
C_4\ve^{\frac{(N-2\alpha)(2-q)+4\alpha}{2} }+C_5\ve^{\frac{q(N-2\alpha)}{2}}&{\rm if}\,\,q\not=\frac{N}{N-2\alpha},\\
C_3\ve^{\frac{(N-2\alpha)(2-q)+4\alpha}{2} }+C_6\ve^{\frac{q(N-2\alpha)}{2}}+C_0\ve^{\frac{q(N-2\alpha)}{2}}|{\rm ln} \ve|&{\rm if}\,\,q=\frac{N}{N-2\alpha}
 \end{array}
 \right.
\end{eqnarray*}
for all $z\in M$. Hence, we have that
\[
\int_\Omega |w_{\varepsilon,z}|^{q}dx=o(\varepsilon^{\frac{N-2\alpha}{2}})
\]
uniformly for $z\in M$, where we have used the fact $1<q<(N+2\alpha)(N-2\alpha)$. $\Box$

Next, we have the following result.

\begin{lemma}\label{l3.3}
Let $ \Lambda_0 $ as defined in Lemma \ref{l3.1}, then, for $ \lambda<\Lambda_0 $,
\[
\sup_{t\ge0}\Phi_\lambda(u_\lambda^++tw_{\varepsilon,z})<d_\lambda:=c_\lambda^++\frac{\alpha}{N}S_\alpha^{\frac{N}{2\alpha}}
\]
uniformly for $ z\in M $.
\end{lemma}

{\bf Proof.} The proof this lemma is very similar to Lemma 3.2 in \cite{CW}. As Lemma 3.2 in \cite{CW}, we first can get the following inequality
\[
\Phi_\lambda(u_\lambda^++tw_{\varepsilon,z})\le \Phi_\lambda(u^+_\lambda)+ J_\lambda(tw_{\ve,z}),
\]
where 
\[
J_\lambda(v)=\frac{1}{2}\|v\|_{X_0^\alpha(\Om)}+C\int_{\Om}v^qdx-\frac{1}{2^*_\alpha}\int_{\Om}b[(u^+_\lambda+v)^{2^*_\alpha}-(u^+_\lambda)^{2^*_\alpha}-2^+_\alpha(u^+_\lambda)^{2^*_\alpha-1}]dx.
\]
Then, by Theorem \ref{t3.1} (i), we just need to prove that
\[
\sup_{t\ge0}J_\lambda(tw_{\ve,z})<\frac{\alpha}{N}S_\alpha^{\frac{N}{2\alpha}}
\]
uniformly for $z\in M$. Applying Lemma \ref{l3.2} and following a similar argument as Lemma 3.2 in \cite{CW}, we can obtain that there exists $t_0>0$ and a sufficiently small $\ve_0$ such that
\[
J_\lambda(tw_{\ve,z})\le 0<\frac{\alpha}{N}S_\alpha^{\frac{N}{2\alpha}}\quad{\rm for\,\,all\,\,}t\in [t_0,\infty),\,\,z\in M\,\,{\rm and}\,\,0<\ve<\ve_0,
\]
and
\[
\max_{t\in [0,t_0]}J_\lambda(tw_{\ve,z})<\frac{\alpha}{N}S_\alpha^{\frac{N}{2\alpha}}.
\]
This completes the proof. $\Box$\\

Next, by using Lemma \ref{l3.3}, we can find a positive solution in  $\mathcal{N}_\lambda^- $ if $ \lambda<\Lambda_0 $.

\begin{teo}\label{t3.2}
Let $ \Lambda_0>0 $ as defined in Lemma \ref{l3.1}, Then, for each $ \lambda<\Lambda_0 $, equation (\ref{1.1}) has a positive solution $ u_\lambda^-\in \mathcal{N}_\lambda^- $.
\end{teo}

{\bf Proof.}  We first show that $c_\lambda^-\le c_\lambda^++\frac{\alpha}{N}S_\alpha^{N/2\alpha}$ and thus we can apply Proposition \ref{p3.1} to obtain a solution. Here we adopt the method
 of Tarantello \cite{T} and Wu \cite{W3}.
By Lemma \ref{l2.5}, we know, for very $ u\in X_0^\alpha( \Omega)\setminus\{0\} $ , that there exists a unique $ t^-=t^-(u)>0 $ such that $ t^-(u)u\in \mathcal{N}_\lambda^- $. So we claim that
\[
\mathcal{N}_\lambda^-=\left\{u\in X_0^\alpha(\Omega)\setminus\{0\}\,\bigg\vert \,\,\frac{1}{\|u\|_{X_0^\alpha(\Omega)}}t^-\left(\frac{u}{\|u\|_{X_0^\alpha(\Omega)}}\right)=1\right\}.
\]
In fact, for $ u\in \mathcal{N}_\lambda^- $, let $w=\frac{u}{\|u\|_{X_0^\alpha(\Omega)}}$. Then there exists a unique $ t^-(w)>0 $ such that $ t^-(w)w\in \mathcal{N}_\lambda^- $ or $ \frac{u}{\|u\|_{X_0^\alpha(\Omega)}}t^-\left(\frac{u}{\|u\|_{X_0^\alpha(\Omega)}}\right)\in \mathcal{N}_\lambda^- $. Since $ u\in \mathcal{N}_\lambda^- $, we have $ \frac{1}{\|u\|_{X_0^\alpha(\Omega)}}t^-\left(\frac{u}{\|u\|_{X_0^\alpha(\Omega)}}\right)=1 $. This implies
\[
\mathcal{N}_\lambda^-\subset\left\{u\in X_0^\alpha(\Omega)\setminus\{0\}\,\bigg\vert \,\,\frac{1}{\|u\|_{X_0^\alpha(\Omega)}}t^-\left(\frac{u}{\|u\|_{X_0^\alpha(\Omega)}}\right)=1\right\}.
\]
Conversely, let $ u\in X_0^\alpha(\Omega)\setminus\{0\} $ such that $ \frac{1}{\|u\|_{X_0^\alpha(\Omega)}}t^-\left(\frac{u}{\|u\|_{X_0^\alpha(\Omega)}}\right)=1 $. Then
\[
u=t^-\left(\frac{u}{\|u\|_{X_0^\alpha(\Omega)}}\right)\frac{u}{\|u\|_{X_0^\alpha(\Omega)}}\in\mathcal{N}_\lambda^-.
\]

Next, we let
\[
A_1=\left\{u\in X_0^\alpha(\Omega)\setminus\{0\}\,\bigg\vert \,\,\frac{1}{\|u\|_{X_0^\alpha(\Omega)}}t^-\left(\frac{u}{\|u\|_{X_0^\alpha(\Omega)}}\right)>1\right\}\cup\{0\};
\]
\[
A_2=\left\{u\in X_0^\alpha(\Omega)\setminus\{0\}\,\bigg\vert \,\,\frac{1}{\|u\|_{X_0^\alpha(\Omega)}}t^-\left(\frac{u}{\|u\|_{X_0^\alpha(\Omega)}}\right)<1\right\}.
\]
Then $ \mathcal{N}_\lambda^- $ disconnects $X_0^\alpha(\Omega) $ in two connected components $ A_1 $ and $ A_2 $. Clearly, $ X_0^\alpha(\Omega)\setminus{N}_\lambda^-=A_1\cup A_2 $ and $ \mathcal{N}_\lambda^+ \subset A_1$. Indeed, for $ u\in  \mathcal{N}_\lambda^+ $, there exist unique $ t^-\left(\frac{u}{\|u\|_{X_0^\alpha(\Omega)}}\right) >0$ and $ t^+\left(\frac{u}{\|u\|_{X_0^\alpha(\Omega)}}\right) >0$ such that
\[
t^+\left(\frac{u}{\|u\|_{X_0^\alpha(\Omega)}}\right)< t_{\rm max}<t^-\left(\frac{u}{\|u\|_{X_0^\alpha(\Omega)}}\right)
\]
and $ t^+\left(\frac{u}{\|u\|_{X_0^\alpha(\Omega)}}\right) \frac{u}{\|u\|_{X_0^\alpha(\Omega)}}\in  \mathcal{N}_\lambda^+ $. Since  $ u\in  \mathcal{N}_\lambda^+ $, we have that
$$ \frac{1}{\|u\|_{X_0^\alpha(\Omega)}}t^+\left(\frac{u}{\|u\|_{X_0^\alpha(\Omega)}}\right)=1 .$$
Therefore,
\[
t^-\left(\frac{u}{\|u\|_{X_0^\alpha(\Omega)}}\right)>t^+\left(\frac{u}{\|u\|_{X_0^\alpha(\Omega)}}\right)=\|u\|_{X_0^\alpha(\Omega)}.
\]
This implies $ \mathcal{N}_\lambda^+ \subset A_1$.

Next, we claim that there exists a $ l_0>0 $ such that $ u_\lambda^++l_0w_{\varepsilon,z}\in A_2 $. Firstly, we find a constant $ C_{19}>0 $ such that $ 0<t^-\left(\frac{ u_\lambda^++lw_{\varepsilon,z}}{\| u_\lambda^++lw_{\varepsilon,z}\|_{X_0^\alpha(\Omega)}}\right)<C_{19} $ for each $ l>0 $. Otherwise, there exists a sequence $ \{l_n\} $ such that $ l_n\rightarrow\infty $ and $t^-\left(\frac{ u_\lambda^++l_nw_{\varepsilon,z}}{\| u_\lambda^++l_nw_{\varepsilon,z}\|_{X_0^\alpha(\Omega)}}\right)\rightarrow\infty$. Let $ v_n=\frac{ u_\lambda^++l_nw_{\varepsilon,z}}{\| u_\lambda^++l_nw_{\varepsilon,z}\|_{X_0^\alpha(\Omega)}} $. Since $ t^-(v_n)v_n\in \mathcal{N}_\lambda^- $, by the Lebesgue dominated convergence theorem,
\begin{eqnarray*}
\int_\Omega b|v_n|^{2_\alpha^*}dx&=&\frac{1}{\| u_\lambda^++l_nw_{\varepsilon,z}\|^{2_\alpha^*}_{X_0^\alpha(\Omega)}}\int_\Omega b(x)| u_\lambda^++l_nw_{\varepsilon,z}|^{2_\alpha^*}dx\\
&=&\frac{1}{\| \frac{u_\lambda^+}{l_n}+w_{\varepsilon,z}\|^{2_\alpha^*}_{X_0^\alpha(\Omega)}}\int_\Omega b(x)\left| \frac{u_\lambda^+}{l_n}+w_{\varepsilon,z}\right|^{2_\alpha^*}dx\\
&\rightarrow&\frac{1}{\| w_{\varepsilon,z}\|^{2_\alpha^*}_{X_0^\alpha(\Omega)}}\int_\Omega \left| w_{\varepsilon,z}\right|^{2_\alpha^*}dx>0
\end{eqnarray*}
as $ n\rightarrow\infty $ and
\begin{eqnarray*}
\Phi_\lambda(t^-(v_n)v_n)&=&\frac{(t^-(v_n))^2}{2}\|v_n\|_{X_0^\alpha(\Omega)}^2-\frac{t^-(v_n))^q}{q}\int_{\Omega}a_\lambda|v_n|^qdx-\frac{t^-(v_n))^{2_\alpha^*}}{2_\alpha^*}\int_\Omega b|v_n|^{2_\alpha^*}dx\\
&\rightarrow &-\infty,
\end{eqnarray*}
as $ n\rightarrow\infty $, which contradicts the fact that $\Phi_\lambda  $ is bounded below on $ {N}_\lambda^- $. Now, we let
\[
l_0=\frac{\left|C_{19}^2-\|u_\lambda^+\|_{X_0^\alpha(\Omega)}^2\right|^{\frac{1}{2}}}{\|w_{\varepsilon,z}\|_{X_0^\alpha(\Omega)}}+1.
\]
Then,
\begin{eqnarray*}
\| u_\lambda^++l_0w_{\varepsilon,z}\|_{X_0^\alpha(\Omega)}^2&=&\| u_\lambda^+\|_{X_0^\alpha(\Omega)}^2+l_0^2\| w_{\varepsilon,z}\|_{X_0^\alpha(\Omega)}^2+2l_0\langle u_\lambda^+,w_{\varepsilon,z}\rangle\\
&>&\| u_\lambda^+\|_{X_0^\alpha(\Omega)}^2+\left|C_{19}^2-\|u_\lambda^+\|_{X_0^\alpha(\Omega)}^2\right|\\
&>& C_{19}^2>\left[t^-\left(\frac{ u_\lambda^++lw_{\varepsilon,z}}{\| u_\lambda^++lw_{\varepsilon,z}\|_{X_0^\alpha(\Omega)}}\right)\right]^2
\end{eqnarray*}
and this implies $ u_\lambda^++l_0w_{\varepsilon,z}\in A_2 $.

Next, we define a path
\[
\gamma(s)=u_\lambda^++sl_0w_{\varepsilon,z}
\]
for $ s\in [0,1] $. Then $ \gamma(0)= u_\lambda^+\in \mathcal{N}_\lambda^+ \subset A_1$ and $ \gamma(1)=u_\lambda^++l_0w_{\varepsilon,z}\in A_2 $. Then there exists a $ s_0\in (0,1) $ such that
\[
\gamma(s_0)=u_\lambda^++s_0l_0w_{\varepsilon,z}\in  \mathcal{N}_\lambda^-.
\]
Therefore, by Lemma \ref{l3.2}, we know that
\[
c_\lambda^-\le \Phi_\lambda(u_\lambda^++s_0l_0w_{\varepsilon,z})<c_\lambda^++\frac{\alpha}{N}S_\alpha^{\frac{N}{2\alpha}}.
\]

Similarly,  by the Ekeland variation principle (see \cite{E}) since  $\Phi_\lambda$ is bounded blow on $\mathcal{N}_\lambda$ as well as on $\mathcal{N}_\lambda^-$, such a minimizing sequence $ \{u_n\}\in \mathcal{N}_\lambda^-$ for $ \Phi_\lambda $ can be established such that
\[
\Phi_\lambda(u_n)=c_\lambda^-+o(1)\quad{\rm and}\quad\Phi_\lambda^\prime(u_n)=o(1)\quad{\rm in}\,\, (X_0^\alpha(\Omega))^*.
\]
By Proposition \ref{p3.1}, there exists a subsequence $ \{u_n\} $ and $ u_\lambda^-\in \mathcal{N}_\lambda^- $ such that $  u_n\rightarrow u_\lambda^-$ strongly in $ X_0^\alpha(\Omega) $, $\Phi_\lambda( u_\lambda^-) =c_\lambda^- $ and $ u_\lambda^- $ is a positive solution of equation (\ref{1.1}) by a similar agurment as in Theorem \ref{t3.1}. $ \Box $\\

{\bf Proof of Theorem \ref{t1.1}}. Together with Theorems \ref{t3.1} and \ref{t3.2}, we obtain Theorem \ref{t1.1}. $ \Box $

\setcounter{equation}{0}
\section{Proof of Theorem \ref{t1.2}}

We first consider the following critical problem
\begin{eqnarray*}
\left\{\begin{array}{l@{\quad }l}
(-\Delta)^\alpha u= |u|^{2^*_\alpha-1}u &{\rm in}\,\,\Omega,\\
u\in X_0^\alpha(\Omega),
 \end{array}
 \right.
 \end{eqnarray*}
and, accordingly, the energy functional $ \Phi^\infty $ in $ X_0^\alpha(\Omega) $ is
\[
 \Phi^\infty(u)=\frac{1}{2}\int_{\R^{2N}}\frac{|u(x)-u(y)|^2}{|x-y|^{N+2\alpha}}dxdy-\frac{1}{2_\alpha^*}\int_\Omega |u|^{2^*_\alpha}dx.
\]
It is easy to check (using the definition of $S_\alpha$) that
\[
\inf_{u\in \mathcal{N}^\infty(\Omega)}\Phi^\infty(u)=\inf_{u\in \mathcal{N}^\infty(\R^N)}\Phi^\infty(u)=\frac{\alpha}{N}S_\alpha^{\frac{N}{2\alpha}},
\]
where
\[
\mathcal{N}^\infty(\R^N)=\{u\in \dot{H}^\alpha(\R^N)\setminus\{0\}\,|\,\,\langle(\Phi^\infty)^\prime(u),u\rangle=0\}
\]
and 
\[
\mathcal{N}^\infty(\Omega)=\{u\in X_0^\alpha(\Omega)\setminus\{0\}\,|\,\,\langle(\Phi^\infty)^\prime(u),u\rangle=0\}
\]
is the Nehari manifold. When $ \lambda=0 $, we write   $ \Phi_\lambda $ (resp. $\mathcal{N}_\lambda  $) as $\Phi_0$ (resp. $\mathcal{N}_0  $). Then, we have the following results.

\begin{lemma}\label{l4.1}
We have that
\[
\inf_{u\in\mathcal{N}_0 }\Phi_0(u)=\inf_{u\in\mathcal{N}^\infty }\Phi^\infty(u)=\frac{\alpha}{N}S_\alpha^{\frac{N}{2\alpha}}.
\]
Moreover, equation (\ref{1.1}) with $ \lambda=0 $ does not admits any positive solution $ u_0 $, for which $ \Phi_0(u_0) =\frac{\alpha}{N}S_\alpha^{\frac{N}{2\alpha}}$.
\end{lemma}

{\bf Proof. } By Lemma \ref{l2.5}, there exists a unique $ t_0(w_{\varepsilon,z}) $ such that
$ t_0(w_{\varepsilon,z})w_{\varepsilon,z}\in \mathcal{N}_0   $ for all $ \varepsilon>0 $, that is,
\beq\label{4.1}
t_0(w_{\varepsilon,z})>t_{{ \rm max}}(w_{\varepsilon,z})=\left(\frac{(2-q)\|w_{\varepsilon,z}\|_{X_0^\alpha(\Omega)}^2}{(2^*_\alpha-q)\int_\Omega b|w_{\varepsilon,z}|^{2^*_\alpha} dx}\right)^{\frac{N-2\alpha}{4\alpha}}
\eqq
and
\beq\label{4.2}
\|t_0(w_{\varepsilon,z})w_{\varepsilon,z}\|_{X_0^\alpha(\Omega)}^2=\int_\Omega a_-|t_0(w_{\varepsilon,z})w_{\varepsilon,z}|^qdx+\int_\Omega b|t_0(w_{\varepsilon,z})w_{\varepsilon,z}|^{2_\alpha^*}dx.
\eqq
Moreover, by Lemma \ref{l3.2} and the boundedness  of $ a_- $, we have
\beq\label{4.3}
\lim_{\varepsilon\rightarrow0}\int_\Omega b|w_{\varepsilon,z}|^{2_\alpha^*}dx=S_\alpha^{\frac{N}{2\alpha}}
\eqq
and
\beq\label{4.4}
\lim_{\varepsilon\rightarrow0}\int_\Omega  a_-|w_{\varepsilon,z}|^qdx=0
\eqq
uniformly for $ z\in M $.

Hence, by (\ref{3.3}) and (\ref{4.1})-(\ref{4.4}), we have
\[
\lim_{\varepsilon\rightarrow0}t_0(w_{\varepsilon,z})=1
\]
uniformly for $ z\in M $. Therefore
\[
\inf_{u\in\mathcal{N}_0 }\Phi_0(u)\le \Phi_0(t_0(w_{\varepsilon,z})w_{\varepsilon,z})\rightarrow \frac{\alpha}{N}S_\alpha^{\frac{N}{2\alpha}}
\]
as $ \varepsilon\rightarrow0 $. So
\[
\inf_{u\in\mathcal{N}_0 }\Phi_0(u)\le\inf_{u\in\mathcal{N}^\infty }\Phi^\infty(u)=\frac{\alpha}{N}S_\alpha^{\frac{N}{2\alpha}}.
\]
Conversely, let $ u\in \mathcal{N}_0 $.  By Lemma \ref{l2.5} and the uniqueness, we have
$\Phi_0(u)=\sup_{t\ge0}\Phi_0(tu)$. Moreover, there exists a unique $ t_u>0 $ such that $ t_uu\in \mathcal{N}^\infty $. Therefore,
\[
\Phi_0(u)\ge \Phi_0(t_uu)\ge \Phi^\infty(t_uu)\ge \frac{\alpha}{N}S_\alpha^{\frac{N}{2\alpha}}.
\]
This implies that $\inf_{u\in\mathcal{N}_0 }\Phi_0(u) \ge \frac{\alpha}{N}S_\alpha^{\frac{N}{2\alpha}}$. Consequently,
\[
\inf_{u\in\mathcal{N}_0 }\Phi_0(u)=\inf_{u\in\mathcal{N}^\infty }\Phi^\infty(u)=\frac{\alpha}{N}S_\alpha^{\frac{N}{2\alpha}}.
\]

Next, we prove that problem (\ref{1.1}) does not admit  any solution $ u_0 $ satisfying $ \Phi_0(u_0) =\inf_{u\in\mathcal{N}_0 }\Phi_0(u)$. We prove it by contradiction. Assume that there exists $ u_0\in \mathcal{N}_0  $ and satisfying $ \Phi_0(u_0) =\inf_{u\in\mathcal{N}_0 }\Phi_0(u)$. As in the proof of Theorem \ref{t3.1}, we may assume $u_0$ is a positive solution.  By Lemma
 \ref{l2.5} gives that $\Phi_0(u_0)=\sup_{t\ge0}\Phi_0(tu_0)$, leading to the conclusion that there must exist a unique $ t_{u_0}>0 $ such that $ t_{u_0}u_0\in \mathcal{N}^\infty $ and
 \begin{eqnarray*}
 \frac{\alpha}{N}S_\alpha^{\frac{N}{2\alpha}}=\inf_{u\in\mathcal{N}_0 }\Phi_0(u)&=&\Phi_0(u_0)\\
 &\ge&\Phi_0(t_ {u_0}u_0)\\
 &\ge&\Phi^\infty(t_ {u_0}u_0)+\frac{1}{2_\alpha^*}\int_\Omega (1-b(x))|t_ {u_0}u_0|^{2_\alpha^*}dx\\
 &>&\frac{\alpha}{N}S_\alpha^{\frac{N}{2\alpha}}+\frac{1}{2_\alpha^*}\int_\Omega (1-b(x))|t_ {u_0}u_0|^{2_\alpha^*}dx.
 \end{eqnarray*}
This implies that
\[
\int_\Omega (1-b(x))|u_0|^{2_\alpha^*}dx<0,
\]
which contradicts  the fact  that $ b\le 1 $ in $ \Omega $. This completes the proof. $ \Box $\\

By using Lemma \ref{l4.1}, we have the following result.

\begin{lemma}\label{l4.2}
Assume that $ \{u_n\} $ is a minimizing sequence for $ \Phi_0 $ in $ \mathcal{N}_0 $. Then,

(1)$\int_\Omega a_-|u_n|^qdx=o(1)$;

(2)$\int_\Omega (1-b)|u_n|^{2_\alpha^*}dx=o(1)$.\\
Moreover, $ \{u_n\}   $ is a $ (PS)_{\frac{\alpha}{N}S_\alpha^{\frac{N}{2\alpha}}} $-sequence for $ \Phi^\infty $ in $ X_0^\alpha(\Omega) $.
\end{lemma}

{\bf Proof.} For each $n$, there is a unique $t_n>0$ such that $t_nu_n\in \mathcal{N}^\infty$, that is,
\[
t_n^2\|u_n\|_{X_0^\alpha(\Omega)}^2=t_n^{2_\alpha^*}\int_\Omega |u_n|^{2_\alpha^*}dx.
\]
By Lemma \ref{l2.5},
\begin{eqnarray*}
\Phi_0(u_n)\ge \Phi_0(t_nu_n)&=&\Phi^\infty(t_nu_n)-\frac{t_n^q}{q}\int_\Omega a_-|u_n|^qdx+\frac{t_n^{2_\alpha^*}}{2_\alpha^*}\int_\Omega (1-b )|u_n|^{2_\alpha^*}dx\\
&\ge&\frac{\alpha}{N}S_\alpha^{\frac{N}{2\alpha}}-\frac{t_n^q}{q}\int_\Omega a_-|u_n|^qdx+\frac{t_n^{2_\alpha^*}}{2_\alpha^*}\int_\Omega (1-b )|u_n|^{2_\alpha^*}dx.
\end{eqnarray*}
Since $\Phi_0(u_n)=\frac{\alpha}{N}S_\alpha^{\frac{N}{2\alpha}}+o(1)$ by Lemma \ref{l3.2}, we have
\[
\frac{t_n^q}{q}\int_\Omega a_-|u_n|^qdx=o(1)
\]
and
\[
\frac{t_n^{2_\alpha^*}}{2_\alpha^*}\int_\Omega (1-b )|u_n|^{2_\alpha^*}dx=o(1).
\]
Next, we prove that there exists $c_0>0$ such that $t_n>c_0$ for all $n$. Suppose the contrary. Then we may assume $t_n\rightarrow0$
 as $n\rightarrow\infty$. Since $\Phi_0(u_n)=\frac{\alpha}{N}S_\alpha^{\frac{N}{2\alpha}}+o(1)$ (see Lemma \ref{4.1}), we know that $\|u_n\|_{X_0^\alpha(\Omega)}$
 is uniformly bounded by Lemma \ref{l2.1}. Hence, $\|t_nu_n\|_{X_0^\alpha(\Omega)}\rightarrow0$ and
 \begin{eqnarray*}
 \Phi^\infty(t_nu_n)&=&\frac{t_n^2}{2}\|u_n\|_{X_0^\alpha(\Omega)}^2-\frac{t_n^{2_\alpha^*}}{2_\alpha^*}\int_\Omega |u_n|^{2_\alpha^*}dx\\
 &=& \frac{\alpha}{N}\|t_nu_n\|_{X_0^\alpha(\Omega)}^2\rightarrow0
 \end{eqnarray*}
 as $n\rightarrow\infty$. This contradict $ \Phi^\infty(t_nu_n)\ge \frac{\alpha}{N}S_\alpha^{\frac{N}{2\alpha}}>0$.

Therefore,
\[
\int_\Omega a_-|u_n|^qdx=o(1)
\]
and
\[
\int_\Omega (1-b )|u_n|^{2_\alpha^*}dx=o(1).
\]
This implies
\[
\|u_n\|_{X_0^\alpha(\Omega)}^2=\int_\Omega |u_n|^{2_\alpha^*}dx+o(1)
\]
and
\[
\Phi^\infty(u_n)=\frac{\alpha}{N}S_\alpha^{\frac{N}{2\alpha}}+o(1).
\]
Then, by a similar argument as Lemma 7 in \cite{WW}, we have $ \{u_n\}   $ is a $ (PS)_{\frac{\alpha}{N}S_\alpha^{\frac{N}{2\alpha}}} $-sequence for $ \Phi^\infty $ in $ X_0^\alpha(\Omega) $. $\Box$\\

Next, for a positive $ d $, we consider the filtration of the Nehari manifold $ \mathcal{N}_0 $ with
\[
\mathcal{N}_0(d)=\left\{u\in \mathcal{N}_0\,\bigg\vert\,\Phi_0(u)\le \frac{\alpha}{N}S_\alpha^{\frac{N}{2\alpha}}+d\right\}
\]
and the function
\[
F(u)=\frac{\int_\Omega x|u|^{2_\alpha^*}dx}{\int_\Omega |u|^{2_\alpha^*}dx}.
\]
With these notations, we have the following result.
\begin{lemma} \label{l4.3}
For each $ 0<\delta<r_0 $, there exists $ d_\delta >0$ such that
\[
F(u)\in M_\delta\quad{\rm for\,\,all}\,\,u\in \mathcal{N}_0(d_\delta).
\]
\end{lemma}

{\bf Proof.} We prove it by contradiction. Assume that there exists a sequence $ \{u_n\}\subset \mathcal{N}_0 $ and $ \delta_0<r_0 $ such that $ \Phi_0(u_n)\le \frac{\alpha}{N}S_\alpha^{\frac{N}{2\alpha}}+o(1) $ and $ F(u_n)\not\in M_{\delta_0} $ for all $ n $.

By Lemma \ref{l4.2}, we know that $ \{u_n\} $ is also a $ (PS)_{\frac{\alpha}{N}S_\alpha^{\frac{N}{2\alpha}}} $-sequence for $ \Phi^\infty $ in $ X_0^\alpha(\Omega) $.
Clearly, $ \|u_n\|_{X_0^\alpha(\Omega)} $ is bounded and thus there exists a subsequence
$ \{u_n\} $ and $ u_0\in X_0^\alpha(\Omega) $ such that $ u_n\rightharpoonup u_0 $ in $ X_0^\alpha(\Omega) $. Since $ \Omega $ is bounded, we have $ u_0\equiv0 $. Therefore, by the concentration-compactness principle (see Theorem 6 in \cite{PP}), there exists two sequences $ \{x_n\}\subset\Omega $ and $ \{R_n\}\subset\R^+ $ with $ x_0\in \bar{\Omega} $ such that $ x_n\rightarrow x_0 $ and $ R_n\rightarrow\infty $ and
\[
\|u_n(x)-R_n^{\frac{N-2\alpha}{2}}u_0(R_n(x-x_n))\|_{L^{2_\alpha^*}(\R^N)}\rightarrow0\quad {\rm as}\,\,n\rightarrow\infty,
\]
where $ u_0 $ is defined as (\ref{2.2}).

Therefore,
\begin{eqnarray*}
F(u_n)&=&\frac{\int_\Omega x|u_n|^{2_\alpha^*}dx}{\int_\Omega |u_n|^{2_\alpha^*}dx}\\
&=&\frac{\int_\Omega x\left| R_n^{\frac{N-2\alpha}{2}}u_0(R_n(x-x_n))\right|^{2_\alpha^*}dx}{\int_\Omega \left| R_n^{\frac{N-2\alpha}{2}}u_0(R_n(x-x_n))\right|^{2_\alpha^*}dx}+o(1)\\
&=&\frac{\int_\Omega \left(\frac{x}{R_n}+x_n\right)| u_0(x)|^{2_\alpha^*}dx}{\int_\Omega | u_0(x)|^{2_\alpha^*}dx}+o(1)\\
&=&x_0+o(1).
\end{eqnarray*}

Next, we show that $ x_0\in  M_{\delta_0}$. Since $ \{u_n\} $ is a minimizing sequence for $ \Phi_0 $ in $ \mathcal{N}_0 $, by Lemma \ref{l4.2},
\begin{eqnarray*}
0=\lim_{n\rightarrow\infty}\int_{\Omega}(1-b)|u_n|^{2_\alpha^*}dx&=&\lim_{n\rightarrow\infty}\int_{\Omega}(1-b)| R_n^{\frac{N-2\alpha}{2}}U_0(R_n(x-x_n))|^{2_\alpha^*}dx\\
&=&\lim_{n\rightarrow\infty}\int_{\Omega}\left(1-b\left(\frac{x}{R_n}+x_n\right)\right)| U_0(x)|^{2_\alpha^*}dx\\
&=&(1-b(x_0))S_\alpha^{\frac{N}{2\alpha}}.
\end{eqnarray*}
This implies that $ b(x_0)=\max_{x\in\bar{\Omega}}b(x)\equiv 1 $ and thus $ x_0\in M $,
which contradicts our assumption. This completes the proof. $ \Box $\\

We now proceed to consider the filtration of the manifold $ \mathcal{N}_\lambda^- $ with
\[
\mathcal{N}_\lambda(c)=\left\{u\in \mathcal{N}_\lambda^-\,|\,\Phi_\lambda(u)\le c\right\}.
\]
We can prove that
\begin{lemma} \label{l4.4}
For each $ 0<\delta<r_0 $, there exists $ 0<\Lambda_\delta\le \Lambda_0$ such that, for $ \lambda< \Lambda_\delta$, we have
\[
F(u)\in M_\delta\quad{\rm for\,\,all}\,\,u\in\mathcal{N}_\lambda(d_\lambda),
\]
where $ d_\lambda $ is defined as in Lemma \ref{l3.3}.
\end{lemma}

{\bf Proof.} For $ u\in\mathcal{N}_\lambda(d_\lambda) $ and thus $ u\in  \mathcal{N}_\lambda^-$, by
(\ref{2.8}) and Lemma \ref{l2.5}, there exists a unique $ t_u>t_{\rm max} (u)$ such that $ t_uu\in \mathcal{N}_0$.  Therefore, by the H\"older and Sobolev inequalities,
\begin{eqnarray}\nonumber
\Phi_\lambda(u)&=&\sup_{t\ge t_{\rm max} (u)}\Phi_\lambda(tu)\\\nonumber
&\ge&\Phi_\lambda(t_uu)\\\nonumber
&=& \Phi_0(t_uu)-\frac{\lambda t_u^q}{q}\int_\Omega a_+|u|^q dx\\\label{4.5}
&\ge& \Phi_0(t_uu)-\frac{\lambda t_u^q}{q}\|a_+\|_{L^{q^*}(\Omega)}S_\alpha^{-\frac{q}{2}}\|u\|_{X_0^\alpha(\Omega)}^q.
\end{eqnarray}
Next, we prove that there exists a positive constant $ \kappa_0 $ independent of $ u $ such that
$ t_u\le \kappa_0 $. In fact, by (\ref{2.8}) and the Sobolev inequality,
\beq\label{4.6}
\|u\|_{X_0^\alpha(\Omega)}^2<\frac{2_\alpha^*-q}{2-q}\int_\Omega b|u|^{2_\alpha^*}dx
\le \frac{2_\alpha^*-q}{2-q}S_\alpha^{-\frac{2_\alpha^*}{2}}\|u\|_{X_0^\alpha(\Omega)}^{2_\alpha^*},
\eqq
and then,
\beq\label{4.7}
\|u\|_{X_0^\alpha(\Omega)}\ge
\left(\frac{(2-q)S_\alpha^{\frac{2_\alpha^*}{2}}}{2_\alpha^*-q}\right)^{\frac{N-2\alpha}{4\alpha}}.
\eqq
Without loss of generality, we may assume that $ t_u\ge1 $. Since
\[
t_u^{2_\alpha^*}\int_\Omega b|u|^{2_\alpha^*}dx=t_u^2\|u\|_{X_0^\alpha(\Omega)}^2-t_u^q\int_\Omega a_-|u|^q dx\le t_u^2\left(\|u\|_{X_0^\alpha(\Omega)}^2+\int_\Omega |a_-||u|^q \right),
\]
we have
\beq\label{4.8}
t_u\le \left( \frac{\|u\|_{X_0^\alpha(\Omega)}^2+\int_\Omega |a_-||u|^q }{\int_\Omega b|u|^{2_\alpha^*}dx}\right)^{\frac{N-2\alpha}{4\alpha}}.
\eqq
Hence, by (\ref{4.6})-(\ref{4.8}) and the H\"older and Sobolev inequalities,
\begin{eqnarray}\nonumber
t_u &\le&\left[\frac{2_\alpha^*-q}{2-q}\left( 1+\frac{\int_\Omega |a_-||u|^q }{\|u\|_{X_0^\alpha(\Omega)}^2}\right)\right]^{\frac{N-2\alpha}{4\alpha}}\\\nonumber
&\le& \left[\frac{2_\alpha^*-q}{2-q}\left( 1+\frac{\|a_-\|_{L^{q^*}(\Omega)}}{S_\alpha^{\frac{q}{2}}\|u\|_{X_0^\alpha(\Omega)}^{2-q}}\right)\right]^{\frac{N-2\alpha}{4\alpha}}\\\nonumber
&\le& \left[\frac{2_\alpha^*-q}{2-q}\left( 1+\|a_-\|_{L^{q^*}(\Omega)}\left(\frac{2_\alpha^*-q}{(2-q)S_\alpha^{\frac{2_\alpha^*-q}{2-q}}}\right)^{\frac{2-q}{2_\alpha^*-q}}\right)\right]^{\frac{N-2\alpha}{4\alpha}}\\\label{4.9}
&=&\kappa_0.
\end{eqnarray}
Substituting (\ref{4.9}) into (\ref{4.5}), we have that
\[
\Phi_0(t_uu)
\le\Phi_\lambda(u)+\frac{\lambda t_u^q}{q}\|a_+\|_{L^{q^*}(\Omega)}S_\alpha^{-\frac{q}{2}}\|u\|_{X_0^\alpha(\Omega)}^q
\le d\lambda+\frac{\lambda \kappa_0^q}{q}\|a_+\|_{L^{q^*}(\Omega)}S_\alpha^{-\frac{q}{2}}\|u\|_{X_0^\alpha(\Omega)}^q.
\]
Since $ \Phi_\lambda(u)\le d_\lambda<\frac{\alpha}{N}S_\alpha^{\frac{N}{2\alpha}} $, by the proof of Lemma \ref{l2.1},, for each $ 0<\lambda<\Lambda_0 $, there exists a positive constant $ c_0 $ independent of $ \lambda $ such that $ \|u\|_{X_0^\alpha(\Omega)}\le c_0 $ for all $ u\in \mathcal{N}_\lambda(d_\lambda) $. Hence,
\[
\Phi_0(t_uu)\le d_\lambda+\frac{\lambda \kappa_0^q}{q}\|a_+\|_{L^{q^*}(\Omega)}S_\alpha^{-\frac{q}{2}}c_0^q.
\]
Let $ d_\delta>0 $ be as in Lemma \ref{l4.3}. Then there exists $ 0<\Lambda_\delta\le \Lambda_0$ such that, for $ \lambda< \Lambda_\delta$,
\[
\Phi_0(t_uu)\le \frac{\alpha}{N}S_\alpha^{\frac{N}{2\alpha}}+d_\delta.
\]
for all $ u\in  \mathcal{N}_\lambda(d_\lambda)$.
By Lemma \ref{l4.3}, we have $ t_uu \in \mathcal{N}_0(d_\delta)$ and
\[
F(u)=\frac{\int_\Omega x|t_uu|^{2_\alpha^*}dx}{\int_\Omega |t_uu|^{2_\alpha^*}dx}=F(t_uu)\in M_\delta
\]
for all $ u\in  \mathcal{N}_\lambda^-(d_\lambda)$. This completes the proof. $ \Box $\\

We recall a multiplicity result for critical points involving Ljusternik-Schnirelman category, which shall apply in proving Theorem \ref{t1.2} (for the proof e.g., see  \cite{MW}).

\begin{teo}\label{t4.1}
Let $ \mathcal{M} $ be a $ C^{1,1} $ complete Riemannian manifold (modelled on a Hilbert space) and assume $ \Psi \in C^1(\mathcal{M},R)$ bounded from below. Let $-\infty<\inf_{\mathcal{M}}\Psi<\sigma<\tau<\infty  $. Suppose that $ \Psi $  satisfies (PS)-condition on the sublevel $ \{u\in \mathcal{M}\,|\,\Psi(u)\le \tau\} $ and $ \sigma $ is not a critical level for $ \Psi $. Then there  exists at least $ {\rm cat}_{\Psi^\sigma}(\Psi^\sigma) $ critical points of $ \Psi $ in $ \Psi^\sigma $, where $ \Psi^\sigma=\{u\in \mathcal{M}\,|\,\Psi(u)\le \sigma\} $.
\end{teo}

\begin{lemma}\label{l4.5}
Let $ X $, $ Y $ and $ Z $ be closed sets with $ Y\subset Z $; let $ h_1\in C(X, Z) $ and $ h_2\in C(Y,X) $. Suppose $ h_1\circ h_2 $ is homotopically equivalent to the identity mapping $ id $ in $ Z $, then $ {\rm cat}_{X}(X)\ge {\rm cat}_{Z}(Y) $.
\end{lemma}

{\bf Proof of Theorem \ref{t1.2}.}
We know that $ \Phi_\lambda $ is $ C^1 $ and $ \mathcal{N}_\lambda^- $ is a $ C^{1,1} $ complete Riemannian manifold. Also $ \Phi_\lambda $ is bounded from below  on $ \mathcal{N}_\lambda^- $ and satisfies $ (PS) $-condition. By Theorem \ref{t4.1}, $ \Phi_\lambda $ has at least  $ {\rm cat}_{\mathcal{N}_\lambda(d_\lambda)}(\mathcal{N}_\lambda(d_\lambda)) $ critical points.

Let $ G_\varepsilon (x) =u_\lambda^++t_zw_{\varepsilon,z}\in \mathcal{N}_\lambda^-$, where $  t_z>0$ depends on $ z $ (see the proof of Theorem \ref{t3.2}). Then,
\[
\Phi_\lambda(u_\lambda^++t_zw_{\varepsilon,z})<d_\lambda=c^+_\lambda+\frac{\alpha}{N}S_\alpha^{\frac{N}{2\alpha}}
\]
for all $ z\in M $ and $ \lambda<\Lambda_0 $, see Lemma \ref{l3.3}. This implies $ u_\lambda^++t_zw_{\varepsilon,z}\in \mathcal{N}_\lambda(d_\lambda) $. Therefore, by Lemma \ref{l4.4}, we know $F\circ G_\varepsilon: M\rightarrow M_\delta$ is well defined. Next, we show that $ F\circ G_\varepsilon $ is hompotopically equivalent to the identity mapping on $ M_\delta $. In fact,
\begin{eqnarray*}
F(G_\varepsilon (z))&=&\frac{\int_\Omega x| G_\varepsilon (z)|^{2_\alpha^*}dx}{\int_\Omega | G_\varepsilon (z)|^{2_\alpha^*}dx}\\
&=&\frac{\int_\Omega x| u_\lambda^++t_zw_{\varepsilon,z}|^{2_\alpha^*}dx}{\int_\Omega | u_\lambda^++t_zw_{\varepsilon,z}|^{2_\alpha^*}dx}\\
&=&\frac{\int_\Omega (x+z)| u_\lambda^+(x+z)+t_z\eta(x)U_\varepsilon(x)|^{2_\alpha^*}dx}{\int_\Omega | u_\lambda^+(x+z)+t_z\eta(x)U_\varepsilon(x)|^{2_\alpha^*}dx}\\
&=&z+\frac{\int_\Omega z| u_\lambda^+(x+z)+t_z\eta(x)U_\varepsilon(x)|^{2_\alpha^*}dx}{\int_\Omega | u_\lambda^+(x+z)+t_z\eta(x)U_\varepsilon(x)|^{2_\alpha^*}dx}\\
&=&z+o(1) \quad {\rm as }\,\,\varepsilon\rightarrow0.
\end{eqnarray*}
Applying Lemma \ref{l4.5} with $ X= \mathcal{N}_\lambda(d_\lambda)$, $ Y=M $, $ Z=M_\delta $,
$ h_1=F $ and $ h_2=G_\varepsilon $, we have $ {\rm cat}_{\mathcal{N}_\lambda(d_\lambda)}(\mathcal{N}_\lambda(d_\lambda)) \ge  {\rm cat}_{M_\delta}(M)$.

Finally, combining the above results with Theorem \ref{t3.1}, we know problem (\ref{1.1}) has at least $ {\rm cat}_{M_\delta}(M) +1$ solutions. $ \Box $

\setcounter{equation}{0}
\section{ Acknowledgements}
The authors would like to express their thanks to the editor and referee for their valuable comments and suggestions.
A. Quaas was partially supported by Fondecyt Grant No. 1151180 Programa Basal, CMM. U. de Chile and Millennium Nucleus Center for Analysis of PDE NC130017.
A. Xia was supported by the "Programa de Iniciaci\'on a la  Investigaci\'on Cient\'ifica" (PIIC) UTFSM 2015.


\begin{thebibliography}{99}

\bibitem{A} D. Applebaum, L\'evy processes and stochastic calculus, Second edition. Cambridge Studies in Advanced Mathematics, 116. Cambridge University Press, Cambridge, 2009.

\bibitem{B1} J. Bertoin, L\'evy processes, Cambridge Tracts in Mathematics, 121. Cambridge University Press, Cambridge, 1996.

\bibitem{B}
 K. Brown and Y. Zhang, The Nehari manifold for a semilinear elliptic equation with a sign-changing weight function. J. Differential Equations 193 (2003), no. 2, 481-499.


 \bibitem{BCSS}
B. Barrios, E. Colorado, R. Servadei and F. Soria, A critical fractional equation with concave-convex power nonlinearities. Ann. Inst. H. Poincar\'e Anal. Non Lin\'eaire 32 (2015), no. 4, 875-900.

 \bibitem{BCPS}
B. Barrios, E. Colorado, A. de Pablo, and U. S\'anchez, On some critical problems for the fractional Laplacian operator. J. Differential Equations 252 (2012), no. 11, 6133-6162.

\bibitem{BL}
H. Brezis and E. Lieb, A relation between pointwise convergence of functions and convergence
of functionals, Proc. Amer. Math. Soc., 88 (1983), no. 3, 486-490.


\bibitem{CDS}
E. Colorado,  A. de Pablo and U. S\'anchez, Perturbations of a critical fractional equation. Pacific J. Math. 271 (2014), no. 1, 65-85.

\bibitem{CW}
 C.Y. Chen and T.F. Wu, Multiple positive solutions for indefinite semilinear elliptic problems involving a critical Sobolev exponent. Proc. Roy. Soc. Edinburgh Sect. A 144 (2014), no. 4, 691-709.

\bibitem{C} D.C. Clark, A variant of Lusternik-Schnirelman theory, Indiana Univ. Math. J. 22 (1972) 65-74.

\bibitem{CT1} R. Cont and P. Tankov, Financial modelling with jump processes, Chapman  and Hall/CRC Financial Mathematics Series. Boca Raton, FI, 2004.

\bibitem{CT}
A. Cotsiolis and N. Tavoularis, Best constants for Sobolev inequalities for higher order fractional derivatives, J. Math. Anal. Appl. 295 (2004), 225-236.



 \bibitem{DP}
P. Dr\'abek and S. Pohozaev, Positive solutions for the p-Laplacian: application of the fibrering method. Proc. Roy. Soc. Edinburgh Sect. A 127 (1997), no. 4, 703-726.

\bibitem{DPV}
E. Di Nezza, G. Palatucci and E. Valdinoci, Hitchhiker's guide to the fractional Sobolev spaces. Bull. Sci. Math. 136 (2012), no. 5, 521-573.

\bibitem{E}
I. Ekeland, On the variational principle, J. Math. Anal. Appl. 17 (1974) 324-353.

\bibitem{J}
I.M. James. On category, in the sense of Ljusternik-Schnirelmann. Topology 17 (1978), 331-348.

\bibitem{MW}
J. Mawhin and M. Willem, Critical Point Theory and Hamiltonian Systems, Springer-
Verlag, Berlin/New York, 1989.

\bibitem{MT}  A. Majda and E. Tabak, A two-dimensional model  for quasigeostrophic flow: comparison with the two-dimensional Euler flow, Nonlinear
phenomena in ocean dynamics (Los Alamos, NM, 1995). Phys. D 98 (1996), no. 2-4, 515522.


\bibitem{N}
Z. Nehari, On a class of nonlinear second-order differential equations, Trans. Amer. Math. Soc. 95 (1960) 101-123.

\bibitem{P}
F.O. de Paiva, Nonnegative solutions of elliptic problems with sublinear indefinite nonlinearity. J. Funct. Anal. 261 (2011), no. 9, 2569-2586.


\bibitem{PP}
G. Palatucci and A. Pisante, Improved Sobolev embeddings, profile decomposition, and
concentration-compactness for fractional Sobolev spaces. Calc. Var. Partial Differential Equations 50 (2014), no. 3-4, 799-829.

\bibitem{P1} S.I. Pohozaev,  An approach to nonlinear equations. (Russian) Dokl. Akad. Nauk SSSR 247 (1979), no. 6, 1327?1331.

\bibitem{R}  P. Rabinowitz, 
Variational methods for nonlinear eigenvalue problems of nonlinear problems, Edizioni Cremonese, Rome, 1974, 139-195. 

\bibitem{SV}
R. Servadei and E. Valdinoci, Mountain pass solutions for non-local elliptic operators. J. Math. Anal. Appl. 389 (2012), no. 2, 887-898.

\bibitem{SV1}
 R. Servadei and  E. Valdinoci, Variational methods for non-local operators of elliptic type. Discrete Contin. Dyn. Syst. 33 (2013), no. 5, 2105-2137.

\bibitem{SV2}
 R. Servadei and  E. Valdinoci,  A Brezis-Nirenberg result for non-local critical equations in low dimension. Commun. Pure Appl. Anal. 12 (2013), no. 6, 2445-2464.

\bibitem{T}
G. Tarantello, On nonhomogeneous elliptic equations involving critical Sobolev exponent, Ann. Inst. H. Poincar\'e Anal. Non Lin\'eaire 9 (1992) 243-261.

\bibitem{V} E. Valdinoci, From the long jump random walk to the fractional Laplacian, Bol. Soc. Esp. Mat. Apl. Se Ma, 49 (2009), 3344.

\bibitem{VIKH} L. Vlahos, H. Isliker, Y. Kominis and K. Hizonidis, Normal and anomalous Diffusion: a tutorial. In Order and chaos, 10th volume, T. Nountis (ed.),
Patras University Press, 2008.

\bibitem{WW}
H.C. Wang and T.F. Wu, Symmetry breaking in a bounded symmetry domain. Nonlinear Differential Equations and Applications NoDEA 11 (2004), 361-377.

\bibitem{W1}
T.F. Wu, Multiple positive solutions for a class of concave-convex elliptic problems in involving sign-changing, J. Funct. Anal. 258 (2010) 99-131.

\bibitem{W2}
T.F. Wu, Three positive solutions for Dirichlet problems involving critical Sobolev exponent and sign-changing weight, J. Differential Equations 249 (2010) 1549-1578.

\bibitem{W3}
T.F. Wu, Multiplicity results for a semi-linear elliptic equation involving sign-changing weight function. Rocky Mountain J. Math. 39 (2009), no. 3, 995-1011.

\bibitem{Y}
X. Yu, The Nehari manifold for elliptic equation involving the square root of the laplacian, J. Differential
Equations 252 (2012), 1283-1308.



\end{thebibliography}
\end{document}